\documentclass[12pt, reqno]{amsart}
\usepackage{pifont}
\usepackage{mathrsfs}
\usepackage{geometry}
\usepackage{titletoc}
\usepackage{stix2}

\usepackage{amsmath}
\usepackage{amssymb} 
\usepackage{enumitem} 
\usepackage{mathtools} 
\usepackage[table]{xcolor} 
\usepackage[all]{xy} 
\usepackage{tikz} 
\usepackage{tikz-cd}
\usepackage{indentfirst} 
\usepackage{babel} 
\usepackage{setspace} 

\usepackage[colorlinks,linkcolor=red,anchorcolor=green,citecolor=blue]{hyperref} 
\hypersetup{linktocpage = true} 

\usepackage{rotating} 

\usepackage{ytableau} 
\usepackage{longtable} 
\newcolumntype{M}[1]{>{\centering\arraybackslash}m{#1}} 

\usepackage{fancyhdr}

\newcommand\bp{{\bar\partial}}

\theoremstyle{plain}
\newtheorem{thm}{Theorem}[section]
\newtheorem{lemma}[thm]{Lemma}
\newtheorem{prop}[thm]{Proposition}
\newtheorem{cor}[thm]{Corollary}
\newtheorem{defn}[thm]{Definition}

\newtheorem{conjecture}[thm]{Conjecture}

\theoremstyle{definition}
\newtheorem{example}[thm]{Example}
\newtheorem{remark}[thm]{Remark}

\newcommand{\btheorem}{\begin{thm}}
    \newcommand{\etheorem}{\end{thm}}
\newcommand{\bproposition}{\begin{prop}}
    \newcommand{\eproposition}{\end{prop}}
\newcommand{\bdefinition}{\begin{defn}}
    \newcommand{\edefinition}{\end{defn}}
\newcommand{\bcorollary}{\begin{cor}}
    \newcommand{\ecorollary}{\end{cor}}
\newcommand{\bproof}{\begin{proof}}
    \newcommand{\eproof}{\end{proof}}
\newcommand{\bremark}{\begin{remark}}
    \newcommand{\eremark}{\end{remark}}
\newcommand{\eexample}{\end{example}}
\newcommand{\bexample}{\begin{example}}

\newcommand{\elemma}{\end{lemma}}
\newcommand{\blemma}{\begin{lemma}}

\newcommand{\p}{\partial}

\renewcommand{\bar}{\overline}
\newcommand{\eps}{\varepsilon}

\renewcommand{\phi}{\varphi}

\newcommand{\beq}{\begin{equation}}
\newcommand{\eeq}{\end{equation}}
\newcommand{\ee}{\end{eqnarray*}}
\newcommand{\be}{\begin{eqnarray*}}

\newcommand{\bd}{\begin{enumerate}}
    \newcommand{\ed}{\end{enumerate}}

\newcommand{\qtq}[1]{\quad\mbox{#1}\quad}
\renewcommand{\bp}{\bar{\partial}}

\newcommand{\ts}{\otimes}

\renewcommand{\S}{{\mathbb S}}

\newcommand{\C}{{\mathbb C}}

\renewcommand{\P}{{\mathbb P}}

\newcommand{\R}{{\mathbb R}}

\newcommand{\om}{\omega}
\newcommand{\smo}{\sqrt{-1}}

\usepackage{harmony}

\setlist[itemize]{leftmargin=*}
\setlist[enumerate]{leftmargin=*}

\numberwithin{equation}{section} 

\makeatletter

\usepackage{fancyhdr}
\title{First eigenvalue estimates on complete K\"ahler manifolds}

\author{Mingwei Wang}
\author{Xiaokui Yang}

\address{Mingwei Wang,  Qiuzhen College, Tsinghua University, Beijing, 100084, China}
\email{wangmw21@mails.tsinghua.edu.cn}

\address{Xiaokui Yang, Department of Mathematics and Yau Mathematical Sciences Center, Tsinghua University, Beijing, 100084, China}
\email{xkyang@mail.tsinghua.edu.cn}

\begin{document}

    \begin{abstract}  Let $ (M,\omega_g) $ be a complete K\"ahler manifold of complex dimension $n$.
     We prove that if the holomorphic sectional curvature satisfies $\mathrm{HSC}  \geq 2 $,
    then the first eigenvalue $\lambda_1$ of the Laplacian on $(M,\omega_g)$ satisfies $$ \lambda_1 \geq \frac{320(n-1)+576}{81(n-1)+144}.$$
    This result is established through a new Bochner-Kodaira type identity specifically developed for holomorphic sectional curvature.
    \end{abstract}

    \maketitle

\section{Introduction}

The spectral analysis of the Laplacian operator on complete Riemannian manifolds
represents a central theme in differential geometry, with deep connections to curvature,
topology, and geometric analysis. Among the most influential results in this
direction are the fundamental theorems of Lichnerowicz (\cite{Lic58})
and Obata (\cite{Oba62}), which establish precise relationships between spectral properties and geometric constraints:
\btheorem[Lichnerowicz-Obata]
Let $(M,g)$ be a complete Riemannian manifold with real dimension $n$ whose Ricci curvature satisfies $\mathrm{Ric}(g)\geq(n-1)g$. Then:
\bd\item The first eigenvalue of the Laplacian  satisfies
$
\lambda_1\geq n.
$
\item Moreover, if $\lambda_1=n$, then $(M,g)$ is isometric to $\left(\S^n,g_{\mathrm{can}}\right)$.
\ed
\etheorem

\noindent This foundational result has inspired numerous extensions
and generalizations in various directions. There are corresponding statements for manifolds with boundary. For 
the Dirichlet problem the result was achieved by Reilly (\cite{Rei77}), and for the 
Neumann problem by Escobar (\cite{Esc88}).  The development of eigenvalue comparisons began with Cheng's work (\cite{Che75}) on first Dirichlet eigenvalues of balls. Building upon Cheng's upper bound for $\lambda_1$, Yau (\cite{Yau75}) pioneered the search for lower bounds expressed through geometric quantities, successfully deriving an estimate incorporating Ricci curvature bounds, diameter, and volume. He further conjectured that the lower bound should depend solely on Ricci curvature and diameter  which was first confirmed by Li (\cite{Li79}) for nonnegative Ricci curvature cases. The general case was subsequently resolved by Li and Yau (\cite{LY80}, \cite{Li82}) through maximum principle techniques. However, their estimate for manifolds with nonnegative Ricci curvature yielded bounds half as large as observed in examples, prompting their conjecture that the true lower bound should be doubled. This final piece was established by Zhong and Yang (\cite{ZY84}). For a comprehensive treatment of these developments, we refer to \cite{LT91} and the references cited therein. In the context of K\"ahler geometry, Futaki  established in
\cite{Fut88} an important eigenvalue estimate that serves as a
K\"ahler analogue to the classical Riemannian results.
\btheorem[Futaki] \label{Futaki}Let $(M,\omega_g)$ be a complete K\"ahler manifold
with complex dimension $n$. Suppose that the Ricci curvature has a
lower bound $\mathrm{Ric}(g)\geq (n+1)g$, then the first eigenvalue
of the Laplacian satisfies $\lambda_1\geq 2(n+1)$. \etheorem

\noindent In the study of spectral geometry on K\"ahler manifolds,
an important distinction from the Riemannian case concerns the
rigidity phenomena. While Obata's theorem establishes a sharp
characterization of the standard sphere when equality is achieved in
the eigenvalue estimate, the corresponding situation in K\"ahler
geometry exhibits fundamentally different behavior.  It is
well-established in the literature that no direct analogue of
Obata's rigidity theorem holds for general K\"ahler manifolds
satisfying the eigenvalue equality case. This distinction arises
from several fundamental aspects including complex structure
constraints and classification complexity. Notably, recent work by Chu, Wang, and Zhang
\cite{CWZ25}  has established important
rigidity results under certain natural geometric conditions.

\btheorem[Chu-Wang-Zhang] Let $(M,\omega_g)$ be a complete K\"ahler manifold of
complex dimension $n$ with Ricci curvature satisfying $\mathrm{Ric}(g)\geq (n+1)g$.
If  the first eigenvalue  of the Laplacian $\lambda_1=2(n+1)$ and either
\bd \item The holomorphic bisectional curvature is positive ($\mathrm{HBSC}>0$); or
\item  The multiplicity of $\lambda_1$ is at least $n^2+3$;
\ed
then $(M,\om_g)$ is isometric to $(\C\P^n,\omega_{\mathrm{FS}})$.

\etheorem

\noindent The aforementioned results can all be interpreted as comparison theorems with corresponding
 rigidity properties. For a more comprehensive treatment of related topics, we refer readers to
 the works of \cite{Pet99}, \cite{LW05}, \cite{Mun09}, \cite{Mun10},
 \cite{TY12}, \cite{Liu14}, \cite{LY18}, \cite{Fuj18}, \cite{NZ18}, \cite{Lot21}, \cite{DSS21}, \cite{Zha22},
  \cite{DS23}, and \cite{XY24+}, along with the references cited therein.\\

Holomorphic sectional curvature is also a fundamental concept in K\"ahler geometry that encodes rich geometric information, sharing important relationships with both Ricci curvature and holomorphic bisectional curvature in the curvature hierarchy of complex manifolds. The classical Schwarz lemma has been extended to holomorphic sectional curvature, yielding powerful Schwarz-type inequalities in both K\"ahler and Hermitian categories  ( e.g. \cite{Yau75}, \cite{CCL79}, \cite{Roy80}, \cite{NZ18}, \cite{CN22}, \cite{Ni21}, \cite{XYY24+}).  Significant geometric consequences arise from curvature conditions: compact K\"ahler manifolds with negative holomorphic sectional curvature are  hyperbolic and possess ample canonical line bundles, making them K\"ahler-Einstein manifolds of general type (\cite{WY16}, \cite{WY16a}, \cite{TY17}, \cite{DT19}). Moreover, compact K\"ahler manifolds with positive holomorphic sectional curvature exhibit strong topological properties--they are projective, rationally connected, and simply connected (\cite{Yan18}). These results are extended to many other important cases and we refer to \cite{YZ19}, \cite{NZ22}, \cite{Mat22}, \cite{Yan24}, \cite{ZZ23+}, \cite{CLZ24+} and the references  therein.\\

This paper establishes new lower bounds for the first eigenvalue in terms of holomorphic sectional curvature.
To the best of our knowledge, these constitute the first lower bounds obtained through this particular approach.

\btheorem
\label{MainTheoremA}
Let $ (M,\omega_g) $ be a complete K\"ahler manifold of complex dimension $n$.
 If  its holomorphic sectional curvature satisfies $\mathrm{HSC}  \geq 2 $, then the first
 eigenvalue $\lambda_1$ of the Laplacian on $(M,\omega_g)$ satisfies \beq   \lambda_1 \geq \frac{320(n-1)+576}{81(n-1)+144}.\eeq 
 In particular,  for any smooth function $\psi \in \C^\infty(M,\R)$ with $\int \psi d\mathrm{vol}_g=0$,
 \beq \int_M| \nabla \psi|^2 d\mathrm{vol}_g\geq  \frac{320(n-1)+576}{81(n-1)+144} \int_M |\psi|^2  d\mathrm{vol}_g. \eeq 
\etheorem

\noindent The proof of Theorem \ref{MainTheoremA} introduces a new Bochner-Kodaira
 formula connecting holomorphic sectional curvature to
 Hodge-theoretic decompositions of eigenfunction differentials. This reveals a deep
  interplay between curvature and spectral geometry, providing new insights into the
  analytic and geometric structure of K\"ahler manifolds.
  \btheorem\label{identity}
  Let $ (M,\omega_g) $ be a compact K\"ahler manifold and  $ f \in C^{\infty}(M,\R) $ be a smooth  function
  satisfying $ \Delta_{\bp} f = \lambda f $ for some $\lambda\in \R$, then the following identity holds:
  \beq 2\lambda \int_M |\phi|^4 \frac{\omega^n}{n!} = \int_M \left( R(V,\bar V,V,\bar V)
  + |\phi|^2 |\p\bar\p f|^2\right) \frac{\omega^n}{n!} + \left\|\omega_1-\lambda f\phi\right\|^2 + \|\omega_1\|^2,\label{key7}\eeq
   where $ \phi = \p f=\omega_g(\cdot, \bar V)$,  
  $\omega_1 = \{ \p_E \phi, \phi \}$ and $E=T^{*1,0} M$.  In particular, one has
  \beq 2\lambda \int_M |V|^4 \frac{\omega^n}{n!} \geq  \int_M  R(V,\bar V,V,\bar V)\frac{\omega^n}{n!}. \label{inequality} \eeq

  \etheorem

  \noindent
   It is well-known that $\left(\C\P^n, \omega_{\mathrm{FS}}\right)$ has $$\mathrm{Ric}( \omega_{\mathrm{FS}})=(n+1) \omega_{\mathrm{FS}}, \quad \mathrm{HSC}\equiv2 \qtq{ and}  \lambda_1=2(n+1).$$ However, the optimal lower bound of
  $\lambda_1$ for holomorphic sectional curvature exhibits qualitatively different
  behavior from the Ricci curvature case, as demonstrated by the following explicit example.

\bexample
Consider the product manifold $(M,\omega_g)$ defined as:
$$
(\mathbb{C}\mathbb{P}^{d_1},a_1\omega_{\mathrm{FS}})\times\cdots\times(\mathbb{C}\mathbb{P}^{d_k},a_k\omega_{\mathrm{FS}})\,,
$$
where $a_i> 0$ and $\sum\limits_{i=1}^k a_i=1$. This manifold has complex dimension $d_1+\cdots +d_k$ and
 holomorphic sectional curvature $\mathrm{HSC}\geq 2$. A straightforward computation shows that the first eigenvalue of the Laplacian on $(M,\omega_g)$ is given by:
\beq 
\lambda_1=2\min_i  \left\{\frac{d_i+1}{a_i}\right\}.\label{producteigenvalue}
\eeq 
\eexample
\noindent This example demonstrates that the optimal lower bound $\Lambda$ of $\lambda_1$ is
independent of the complex dimension $n$ and satisfies  \beq \frac{320}{81}< \Lambda\leq 4.\eeq  We propose:
\begin{conjecture} Let $ (M,\omega_g) $ be an $n$-dimensional compact K\"ahler manifold with holomorphic sectional
curvature $\mathrm{HSC}  \geq 2 $.  Then the first
eigenvalue $\lambda_1\geq 4$.
Moreover, $\lambda_1=4$ if and only if $(M,\omega_g)$ is isometric to $(\C\P^1, \omega_{\mathrm{FS}})$.
\end{conjecture}

\noindent\textbf{Acknowledgements}.  The authors would like to express their gratitude to Zhiyao Xiong and Jiaxuan Fan for discussions. The second author is particularly grateful to Bing-Long Chen and Valentino Tosatti for their insightful comments and suggestions. 

\vskip 2\baselineskip

\section{Proof of the main theorem}\label{2}

Let $ (M,\omega_g) $ be a compact K\"ahler manifold and $ E = T^{*1,0}M $ be  the  holomorphic tangent bundle.
 Given a smooth real function $ f \in C^{\infty}(M,\mathbb{R}) $,
 we define  $ \phi := \p f\in \Gamma(M,E) $ as its $(1,0)$-gradient and let $ V $ be  the dual vector field of $ \bar\phi $ with respect to the K\"ahler metric.
For notational convenience, we introduce the following differential $1$-forms:
\beq \omega_1 = \{ \p_E \phi, \phi \}, \qquad \omega_2 = \{ \bar\p_E \phi, \phi \}, \eeq
where $\p_E$ and $\bp_E$ denote the $(1,0)$ and $(0,1)$ components of the Chern connection on
$E$ respectively. These forms satisfy the  identity
\begin{equation}
\partial|\phi|^2 = \omega_1 + \overline{\omega}_2.
\end{equation}

\noindent We present two computational lemmas.

\blemma
\label{omegatwo} One has  $\omega_2 = I_V \p\bar\p f$ and
\beq \label{omegatwoA} \langle \omega_2, \omega_2 \rangle = - \langle \{\bar\p_E\phi,\bar\p_E\phi\}_E, \phi\wedge\bar\phi\rangle,
 \quad \langle \omega_1, \omega_1 \rangle =  \langle \{\p_E\phi,\p_E\phi\}_E, \phi\wedge\bar\phi\rangle. \eeq
\elemma

\bproof If we write  $ \p\bar\p f = f_{i\bar j}dz^i \wedge d\bar z{}^j $ and $\phi=\phi_i dz^i$, then $\bp_E \phi=f_{i\bar j} d\bar z^j\ts dz^i$ and
\beq \omega_2 = \{f_{i\bar j}d\bar z^j\ts dz^i, \phi_k dz^k \}_E = g^{i\bar k}f_{i\bar j}\bar \phi_k d\bar z{}^j. \eeq
On the other hand,  one has $V=g^{i\bar k}\bar \phi_k\frac{\p}{\p z^i}$ and
\beq I_V \p\bar\p f = g^{i\bar k}\bar \phi_k I_i(f_{\ell\bar j}dz^\ell\wedge d\bar z{}^j) = g^{i\bar k}f_{i\bar j}\bar \phi_k d\bar z{}^j, \eeq
and we establish the  identity $\omega_2 = I_V \p\bar\p f$. For the second identity, one has
\beq \langle \omega_2, \omega_2 \rangle = g^{\ell\bar j}g^{i\bar a}g^{b\bar k}f_{i\bar j}\bar f_{k\bar\ell}\bar \phi_a \phi_b. \label{2.1}\eeq
On the other hand,  one can see clearly that \beq  \{ \bar\p_E \phi,\bar\p_E \phi \}_E = g^{i\bar k}f_{i\bar j}\bar f_{k\bar\ell} d\bar z{}^j \wedge dz^\ell,\eeq
and so
\beq - \langle \{\bar\p_E\phi,\bar\p_E\phi\}_E, \phi\wedge\bar\phi\rangle = g^{\ell\bar a}g^{b\bar j}g^{i\bar k}f_{i\bar j}\bar f_{k\bar\ell}\bar\phi_a\phi_b. \label{2.2}\eeq
Equations \eqref{2.1} and \eqref{2.2} give
\beq \langle \omega_2, \omega_2 \rangle = - \langle \{\bar\p_E\phi,\bar\p_E\phi\}_E, \phi\wedge\bar\phi\rangle.\eeq
For the third identity, if we write $ \p_E \p f = t_{ij}dz^i \otimes dz^j $  where $ t_{ij} = t_{ji} = \frac{\p^2 f}{\p z^i\p z^j}- \Gamma_{ij}^k\frac{\p f}{\p z^k} $,
then by similar computations, we have
\beq \omega_1 = \{t_{ij}dz^i \otimes dz^j, \phi_k dz^k \}_E = g^{i\bar k}t_{ij}\bar \phi_k dz^j, \eeq
and
\beq \langle \omega_1, \omega_1 \rangle = g^{j\bar\ell}g^{i\bar a}g^{b\bar k}t_{ij}\bar t_{k\ell}\bar \phi_a \phi_b=
\langle \{\p_E\phi,\p_E\phi\}_E, \phi\wedge\bar\phi\rangle. \eeq
This completes the proof.\eproof

\blemma\label{dstarcompute} Let $ \phi = \p f \in \Omega^{1,0}(M) $.  For any smooth function $ p \in C^{\infty}(M,\mathbb{R}) $, one has
\beq \p^*\left(p\phi\wedge\bar\phi\right) =  p(\Delta_{\bar\p}f)\bar\phi - p\omega_1 - \langle \phi ,\p p \rangle \bar\phi, \eeq
and
\beq \bar\p{}^*\left(p\phi\wedge\bar\phi\right) =  -p(\Delta_{\bar\p}f)\bar\phi + p\omega_1 + \langle \bar\phi ,\bar\p p \rangle \phi. \eeq
In particular, for any $ \eps > 0 $, one has
\beq \p^*\left(\frac{\phi\wedge\bar\phi}{|\phi|^2 + \eps}\right) = \frac{(\Delta_{\bar\p}f)\bar\phi - \bar\omega_1}{|\phi|^2 + \eps} + \frac{\langle \phi, \bar\omega_2 + \omega_1 \rangle}{(|\phi| + \eps)^2} \cdot\bar\phi, \eeq
and
\beq \bar\p{}^*\left(\frac{\phi\wedge\bar\phi}{|\phi|^2+\eps}\right) =  \frac{-(\Delta_{\bar\p}f)\phi + \omega_1}{|\phi|^2+\eps} - \frac{\langle \bar\phi, \omega_2 + \bar\omega_1 \rangle}{(|\phi|+\eps)^2} \cdot\phi. \eeq
\elemma
\bproof By using the elementary identity $ \p^* = \smo[\Lambda,\bar\p] $, one has
\be  \p^*\left(p\phi\wedge\bar\phi\right) &=& \smo[\Lambda,\bar\p]\left(p\phi\wedge\bar\phi\right) \\&=& \smo\Lambda\bar\p\left(p\phi\wedge\bar\phi\right) - \bar\p(p|\phi|^2) \\
& = & \smo\Lambda\left(p\bar \p\phi\wedge\bar\phi + \bar\p p \wedge \phi \wedge \bar \phi \right) - |\phi|^2 \bar\p p - p\bar\p|\phi|^2. \ee
On the other hand, one has the contraction formula $ \smo\Lambda = g^{i\bar j}I_{\bar j} I_i $ and so
\beq \smo\Lambda (\bar\p\phi \wedge \bar\phi) = g^{i\bar j}I_{\bar j}(I_i \bar\p\phi \wedge \bar \phi) = (\Delta_{\bp} f)\bar \phi + I_V \p\bar\p f. \eeq
Moreover,
\beq \smo\Lambda(\bar\p p \wedge \phi \wedge \bar \phi) = -g^{i\bar j}\phi_i I_{\bar j}(\bar\p p \wedge \bar \phi) = |\phi|^2 \bar\p p - (\bar Vp)\bar \phi. \eeq
Therefore,
\beq \p^*\left(p\phi\wedge\bar\phi\right) = p(\Delta_{\bar\p}f)\bar\phi + pI_V\p\bar\p f -(\bar V p)\bar\phi - p\bar\p|\phi|^2. \eeq
By Lemma \ref{omegatwo}, we obtain
\beq \p^*\left(p\phi\wedge\bar\phi\right) = p(\Delta_{\bar\p}f)\bar\phi - p\omega_1 - \langle \phi ,\p p \rangle \bar\phi,  \eeq
and the second equation can be proved in similar ways.
\eproof

\noindent We now present the proof of Theorem \ref{MainTheoremA}, whose foundational approach builds upon the key ideas developed in the proof of Theorem \ref{identity}.

\bproof[Proof of  Theorem \ref{MainTheoremA}] It is well-known (e.g.
\cite{Tsu57}, \cite{XY24+}) that a complete K\"ahler manifold with
$\mathrm{HSC}\geq 2$ is compact. Let $ \lambda_1 > 0 $ be the first
eigenvalue of  $(M,\omega_g)$ and $ f \in C^{\infty}(M,\mathbb{R}) $
be a nonzero eigenfunction, i.e. $\Delta_d f=(dd^*+d^*d)f=\lambda_1
f$. In particular, we have $\Delta_{\bp } f=\lambda_1 f/2$. We set $\lambda=\lambda_1/2$.  When $n=1$ and assume $n\geq 2$. 

 For any $\eps>0$,   we consider the $(1,1)$-form $ \alpha \in \Omega^{1,1}(M) $
\beq \alpha = \frac{\p\{\bar\p_E\phi,\phi\} + \bar\p\{\p_E\phi,\phi\}}{|\phi|^2+\eps}.\eeq
It is clear that  \beq\alpha = \frac{\{R^E\phi,\phi\} - \{\p_E\phi,\p_E\phi\} - \{\bar\p_E\phi,\bar\p_E\phi\}}{|\phi|^2+\eps}. \eeq
Moreover, if we write $ \phi = \phi_i dz^i $,
\beq R^E\phi = -g^{d\bar\ell}R_{i\bar jk\bar\ell} \phi_d dz^k \otimes (dz^i \wedge d\bar z{}^j). \eeq
Hence, we obtain an identity for $(1,1)$ forms
\beq \{ R^E\phi, \phi \} = -g^{k\bar c}g^{d\bar\ell}R_{i\bar jk\bar\ell} \bar\phi_c\phi_d dz^i \wedge d\bar z{}^j. \eeq
Since $ \phi \wedge \bar \phi = \phi_a\bar\phi_b dz^a \wedge d\bar z{}^b $,
\beq \langle \{ R^E\phi,\phi \}, \phi \wedge \bar\phi \rangle = -g^{i\bar a}g^{b\bar j}g^{k\bar c}g^{d\bar\ell}R_{i\bar jk\bar\ell} \bar\phi_a\phi_b\bar\phi_c\phi_d = -R(V,\bar V,V,\bar V). \eeq
On the other hand, since $ \smo \alpha $ and $ \smo \phi \wedge \bar\phi $ are both real forms,
\be (\phi\wedge \bar\phi,\alpha) & = & (\alpha,\phi\wedge \bar\phi) = \left(\frac{\{R^E\phi,\phi\} - \{\p_E\phi,\p_E\phi\} - \{\bar\p_E\phi,\bar\p_E\phi\}}{|\phi|^2+\eps}, \phi \wedge \bar\phi \right) \\
& = & \int_M -\frac{R(V,\bar V,V,\bar V)}{|\phi|^2+\eps} \frac{\omega^n}{n!} - \left( \{\p_E\phi,\p_E\phi\}, \frac{\phi\wedge \bar\phi}{|\phi|^2+\eps}\right) - \left( \{\bar\p_E\phi,\bar\p_E\phi\}, \frac{\phi\wedge \bar\phi}{|\phi|^2+\eps}\right) \\
& = & -\int_M \frac{R(V,\bar V,V,\bar V)}{|\phi|^2+\eps} \cdot \frac{\omega^n}{n!} - \int_M \frac{|\omega_1|^2}{|\phi|^2+\eps} \cdot \frac{\omega^n}{n!} + \int_{M} \frac{|\omega_2|^2}{|\phi|^2+\eps} \cdot \frac{\omega^n}{n!}, \ee
where the third identity follows from Lemma \ref{omegatwo}. On the other hand, by Lemma \ref{dstarcompute},
\be (\phi\wedge \bar\phi,\alpha) & = & \left( \frac{\phi \wedge \bar\phi}{|\phi|^2+\eps}, \p\omega_2 + \bar\p\omega_1 \right) = \left( \p^* \left(\frac{\phi \wedge \bar\phi}{|\phi|^2+\eps}\right), \omega_2\right) + \left( \bar\p{}^* \left(\frac{\phi \wedge \bar\phi}{|\phi|^2+\eps}\right), \omega_1\right) \\
& = & \int_M \left(\frac{(\Delta_{\bar\p}f)\langle \bar\phi,\omega_2 \rangle - \langle\bar\omega_1,\omega_2 \rangle}{|\phi|^2+\eps} + \frac{\langle \phi, \bar\omega_2 + \omega_1 \rangle \langle\bar\phi, \omega_2\rangle}{\left(|\phi|^2+\eps\right)^2}\right) \frac{\omega^n}{n!} \\
& + & \int_M \left(\frac{-(\Delta_{\bar\p}f)\langle \phi,\omega_1 \rangle + \langle \omega_1,\omega_1 \rangle}{|\phi|^2+\eps} - \frac{\langle \bar\phi, \omega_2 + \bar\omega_1 \rangle\langle\phi, \omega_1\rangle}{\left(|\phi|^2+\eps\right)^2}\right) \frac{\omega^n}{n!} \\
& = & \int_M \left(\frac{(\Delta_{\bar\p}f)\langle \bar\phi,\omega_2 \rangle - \langle\bar\omega_1,\omega_2 \rangle}{|\phi|^2+\eps} + \frac{|\langle \omega_2 , \bar\phi \rangle|^2}{\left(|\phi|^2+\eps\right)^2}\right)\frac{\omega^n}{n!} \\
& + & \int_M \left(\frac{-(\Delta_{\bar\p}f)\langle \phi,\omega_1 \rangle + \langle \omega_1,\omega_1 \rangle}{|\phi|^2+\eps} - \frac{|\langle \omega_1, \phi \rangle|^2}{\left(|\phi|^2+\eps\right)^2}\right) \frac{\omega^n}{n!} \ee
Therefore, one has
\be
0 & = & \int_M \frac{R(V,\bar V,V,\bar V)}{|\phi|^2+\eps} \frac{\omega^n}{n!} + \int_M \left(\frac{(\Delta_{\bar\p}f)\langle \bar\phi,\omega_2 \rangle }{|\phi|^2+\eps} + \frac{|\langle \omega_2 , \bar\phi \rangle|^2}{\left(|\phi|^2+\eps\right)^2} - \frac{\langle \omega_2, \omega_2 \rangle}{|\phi|^2+\eps}\right) \frac{\omega^n}{n!} \\
&& +  \int_M \frac{-(\Delta_{\bar\p}f)\langle \phi,\omega_1 \rangle- \langle \bar\omega_2,\omega_1 \rangle + 2\langle \omega_1,\omega_1 \rangle}{|\phi|^2+\eps} - \frac{|\langle \omega_1, \phi \rangle|^2}{\left(|\phi|^2+\eps\right)^2}  \frac{\omega^n}{n!}. \ee
Let $M_0=\{x\in M\ |\ \phi(x)=0 \}$. It is obvious that the norms of $|V| , |\omega_1|, |\omega_2|$ are uniformly bounded by $C|\phi|$. By using Lebesgue's dominated convergence theorem,
\begin{eqnarray}
0 \nonumber& = & \int_{M\setminus M_0} \frac{R(V,\bar V,V,\bar V)}{|\phi|^2} \frac{\omega^n}{n!} + \int_{M\setminus M_0} \left(\frac{(\Delta_{\bar\p}f)\langle \bar\phi,\omega_2 \rangle }{|\phi|^2} + \frac{|\langle \omega_2 , \bar\phi \rangle|^2}{|\phi|^4} - \frac{\langle \omega_2, \omega_2 \rangle}{|\phi|^2}\right) \frac{\omega^n}{n!} \\
&  &+ \int_{M\setminus M_0} \frac{-(\Delta_{\bar\p}f)\langle \phi,\omega_1 \rangle- \langle \bar\omega_2,\omega_1 \rangle + 2\langle \omega_1,\omega_1 \rangle}{|\phi|^2} - \frac{|\langle \omega_1, \phi \rangle|^2}{|\phi|^4}  \frac{\omega^n}{n!}. \end{eqnarray}
By taking the conjugate, we obtain
\begin{eqnarray}
0 \nonumber& = & \int_{M\setminus M_0} \frac{R(V,\bar V,V,\bar V)}{|\phi|^2} \frac{\omega^n}{n!}\\ &+& \int_{M\setminus M_0} \left(\frac{(\Delta_{\bar\p}f)\left(\langle \bar\phi,\omega_2 \rangle+\langle\omega_2, \bar\phi \rangle\right) }{2|\phi|^2} + \frac{|\langle \omega_2 , \bar\phi \rangle|^2}{|\phi|^4} - \frac{\langle \omega_2, \omega_2 \rangle}{|\phi|^2}\right) \frac{\omega^n}{n!} \label{keyidentity}\\
\nonumber& + & \int_{M\setminus M_0}\left( \frac{-(\Delta_{\bar\p}f)\left(\langle \phi,\omega_1 \rangle+\langle \omega_1,\phi \rangle\right)- \langle \bar\omega_2,\omega_1 \rangle - \langle \omega_1,\bar\omega_2 \rangle+ 4\langle \omega_1,\omega_1 \rangle}{2|\phi|^2} - \frac{|\langle \omega_1, \phi \rangle|^2}{|\phi|^4}\right)  \frac{\omega^n}{n!}. \end{eqnarray}

\noindent
We shall perform computations on $M\setminus M_0$. We set:
\beq \omega_1^{\perp} = \omega_1 - \frac{\langle \omega_1,\phi \rangle}{|\phi|^2}\cdot \phi, \quad \bar\omega{}_2^{\perp} = \bar\omega_2 - \frac{\langle \bar\omega_2,\phi \rangle}{|\phi|^2}\cdot \phi. \eeq
 The third line in the identity \eqref{keyidentity} is replaced via the following combinatorial formula which can be verified in a straightforward way:
\be 0 & \leq & 2\left|\omega_1^\perp - \frac{1}{4}\bar\omega{}_2^\perp \right|^2 + \frac{1}{4|\phi|^2}\left|\langle 2\omega_1 - \Delta_{\bar\p}f\cdot\phi - \bar\omega_2,\phi \rangle \right|^2 \\
& = & |\omega_1^\perp|^2 - \frac{1}{8}|\bar\omega{}_2^\perp|^2 + \frac{1}{4}|2\omega_1^\perp - \bar\omega{}_2^\perp|^2 +  \frac{1}{4|\phi|^2}\left|\langle 2\omega_1 - \Delta_{\bar\p}f\cdot\phi - \bar\omega_2,\phi \rangle \right|^2. \\
& = & |\omega_1^\perp|^2 - \frac{1}{8}|\bar\omega_2^\perp|^2 + \frac{1}{4}\left| 2\omega_1 -\Delta_{\bar\p} f\cdot \phi -\bar\omega_2 \right|^2 \\
& = & - \frac{1}{2}\Delta_{\bar \p}f \left(\langle \omega_1,\phi \rangle + \langle \phi, \omega_1 \rangle \right) - \frac{1}{2} \left(\langle \omega_1,\bar\omega_2 \rangle + \langle \bar\omega_2, \omega_1 \rangle \right) + 2|\omega_1|^2 - \frac{|\langle \omega_1, \phi \rangle|^2}{|\phi|^2} \\
&& + \frac{1}{4}(\Delta_{\bar\p}f)^2|\phi|^2 + \frac{1}{4}\Delta_{\bar \p}f \left(\langle \bar\omega_2,\phi \rangle + \langle \phi, \bar\omega_2 \rangle \right) + \frac{1}{8}|\omega_2|^2 + \frac{|\langle \phi,\bar\omega_2\rangle|^2}{8|\phi|^2}. \ee
By using this inequality,  identity \eqref{keyidentity} is reduced to
\begin{eqnarray} 0 &\geq& \int_{M\setminus M_0} \frac{R(V,\bar V,V,\bar V)}{|\phi|^2} \frac{\omega^n}{n!} \label{key4}\\ \nonumber&+& \int_{M\setminus M_0} \left(\frac{\Delta_{\bar \p}f \left(\langle \bar\omega_2,\phi \rangle + \langle \phi, \bar\omega_2 \rangle \right)}{4|\phi|^2}+  \frac{7|\langle \phi, \bar\omega_2\rangle|^2}{8|\phi|^4}  - \frac{9|\omega_2|^2}{8|\phi|^2} -\frac{1}{4}(\Delta_{\bar\p} f)^2 \right) \frac{\omega^n}{n!}. \end{eqnarray}
We claim that the following inequality holds on $M\setminus M_0$:
\begin{eqnarray} \nonumber&&\frac{\Delta_{\bar \p}f \left(\langle \bar\omega_2,\phi \rangle + \langle \phi, \bar\omega_2 \rangle \right)}{4|\phi|^2}+  \frac{7|\langle \phi, \bar\omega_2\rangle|^2}{8|\phi|^4}  - \frac{9|\omega_2|^2}{8|\phi|^2}\\&\geq& -\frac{9}{16}|\p\bp f|^2-\frac{16(n-1)+27}{80(n-1)+144}(\Delta_{\bp} f)^2.\label{claim}\end{eqnarray} 
Given this inequality, we have
\beq 0 \geq \int_{M\setminus M_0} \frac{R(V,\bar V,V,\bar V)}{|\phi|^2} \frac{\omega^n}{n!}-\int_{M\setminus M_0} \left(\frac{9}{16}|\p\bp f|^2 +\left(\frac{1}{4}+\frac{16(n-1)+27}{80(n-1)+144}\right)(\Delta_{\bp} f)^2\right)\frac{\omega^n}{n!}.\label{key6} \eeq 
Since  $\phi|_{M_0}=0$,   it is well-known that  (e.g. \cite[p.~310]{Eva2010}),  \beq \p\bp f=-\bp\phi=0, \qtq{a.e. on}M_0. \eeq
In particular, one has
\beq \int_{M\setminus M_0} (\Delta_{\bar\p} f)^2\frac{\omega^n}{n!} =\int_{M} (\Delta_{\bar\p} f)^2\frac{\omega^n}{n!} \qtq{and} \int_{M\setminus M_0} |\p\bar\p f|^2\frac{\omega^n}{n!}=\int_{M} |\p\bar\p f|^2\frac{\omega^n}{n!}.  \eeq
Note also that $\Delta_{\bp } f=\lambda f$ and
\beq \int_{M} (\Delta_{\bar\p} f)^2\frac{\omega^n}{n!} = \lambda\int_M f\Delta_{\bar\p} f \frac{\omega^n}{n!} =\lambda \int_M |\bp f|^2 \frac{\omega^n}{n!}=\lambda \int_M |\phi|^2 \frac{\omega^n}{n!},\eeq
and
\beq \int_{M} |\p\bar\p f|^2\frac{\omega^n}{n!}=(\p\bp f, \p\bp f)=(\Delta_{\p} (\bp f), \bp f)=\lambda \int_M |\phi|^2 \frac{\omega^n}{n!}.  \eeq
Finally, we obtain from \eqref{key6} that
\beq 0\geq  2\int_M |\phi|^2 \frac{\omega^n}{n!} -\left(\frac{9}{16}+\frac{1}{4}+\frac{16(n-1)+27}{80(n-1)+144}\right)\lambda \int_M |\phi|^2 \frac{\omega^n}{n!}.\eeq
Therefore, 
\beq \lambda_1=2\lambda \geq \frac{320(n-1)+576}{81(n-1)+144}.\eeq
This completes the proof of Theorem \ref{MainTheoremA}.\\

 Now we prove the  inequality  \eqref{claim}. At a given  point $ q \in M\setminus M_0 $, there is an  orthonormal basis $ \{ e_i \}_{i=1}^n $ of $ T_qM $ such that $ e_n = \frac{V}{|V|} $. Let $ \{ e^i \} $ be the dual frame of $ \{ e_i \} $.   Suppose that  $ A=(a_{i\bar j}) $ is the Hermitian matrix with $ a_{i\bar j} = (\p\bp f)(e_i,\bar e_j) $.  It is easy to see that at point $q$,  $$ \p\bar\p f = \sum_{i, j}a_{i\bar j}e^i\wedge \bar e{}^j , \quad  \omega_2 = I_V \p\bar\p f = |\phi|\sum_i a_{n\bar i} \cdot \bar e{}^i $$ and so
\beq  |\p\bar\p f|^2 = \sum_{i,j=1}^n |a_{i\bar j}|^2, \quad \frac{|\omega_2|^2}{|\phi|^2}= \sum_{i=1}^n |a_{n\bar i}|^2, \quad \frac{|\langle \phi,\bar\omega_2 \rangle|^2}{|\phi|^4}  =\frac{|\bar\omega_2 (V)|^2}{|\phi|^4} = |a_{n\bar n}|^2. \eeq
Therefore, we  obtain
\be &&\frac{\Delta_{\bar \p}f \left(\langle \bar\omega_2,\phi \rangle + \langle \phi, \bar\omega_2 \rangle \right)}{4|\phi|^2}+  \frac{7|\langle \phi, \bar\omega_2\rangle|^2}{8|\phi|^4}  - \frac{9|\omega_2|^2}{8|\phi|^2}\\ &=& -\frac{9}{8}\sum_{i=1}^n |a_{n\bar i}|^2 + \frac{7}{8}|a_{n\bar n}|^2 - \frac{1}{2}\sum_{i=1}^n a_{i\bar i}a_{n\bar n}. \ee
Since one has 
\beq -2\sum_{i=1}^n|a_{n\bar i}|^2 \geq - |a_{n\bar n}|^2 + \sum_{i\neq n}|a_{i\bar i}|^2 - \sum_{i,j=1}^n |a_{i\bar j}|^2. \eeq
Hence,
\be 
&& -\frac{9}{8}\sum_{i=1}^n |a_{n\bar i}|^2 + \frac{7}{8}|a_{n\bar n}|^2 - \frac{1}{2}\sum_{i=1}^n a_{i\bar i}a_{n\bar n} \\
& \geq & \frac{5}{16}|a_{n\bar n}|^2 + \frac{9}{16}\sum_{i\neq n}|a_{i\bar i}|^2 - \frac{1}{2}\sum_{i=1}^n a_{i\bar i}a_{n\bar n} - \frac{9}{16}\sum_{i,j=1}^n |a_{i\bar j}|^2 \\
&= &-\frac{3}{16}|a_{n\bar n}|^2 + \frac{9}{16}\sum_{i\neq n}|a_{i\bar i}|^2  - \frac{1}{2}\left(\sum_{i\neq n}a_{i\bar i}\right)a_{n\bar n} - \frac{9}{16}\sum_{i,j=1}^n |a_{i\bar j}|^2 \\
& \geq & -\frac{3}{16}|a_{n\bar n}|^2 - \frac{1}{2}\left(\sum_{i\neq n}a_{i\bar i}\right)a_{n\bar n}+ \frac{9}{16(n-1)} \left( \sum_{i\neq n} a_{i\bar i}\right)^2 - \frac{9}{16}\sum_{i,j=1}^n |a_{i\bar j}|^2. \ee
For any $\kappa_0\in \R$, one has
\be&&   -\frac{3}{16}|a_{n\bar n}|^2 - \frac{1}{2}\left(\sum_{i\neq n}a_{i\bar i}\right)a_{n\bar n}+ \frac{9}{16(n-1)} \left( \sum_{i\neq n} a_{i\bar i}\right)^2\\
& =& \left(-\frac{3}{16}+\kappa_0\right)|a_{n\bar n}|^2 + \left(-\frac{1}{2}+2\kappa_0\right)\left(\sum_{i\neq n}a_{i\bar i}\right)a_{n\bar n}+ \left(\frac{9}{16(n-1)} + \kappa_0\right) \left( \sum_{i\neq n} a_{i\bar i}\right)^2\\
&&  - \kappa_0\left(\sum_{i=1}^n a_{i\bar i}\right)^2. \ee
Let 
\beq \kappa_0 = \frac{16(n-1)+27}{80(n-1)+144},  \eeq
 be the solution of the discriminant  of the quadratic equation
\beq \Delta=\left(-\frac{1}{2} +2 \kappa_0\right)^2- 4\cdot  \left(-\frac{3}{16} + \kappa_0\right)\left(\frac{9}{16(n-1)} + \kappa_0\right) =0.  \eeq
One can deduce that
\beq \left(-\frac{3}{16}+\kappa_0\right)|a_{n\bar n}|^2 + \left(-\frac{1}{2}+2\kappa_0\right)\left(\sum_{i\neq n}a_{i\bar i}\right)a_{n\bar n} + \left(\frac{9}{16(n-1)} + \kappa_0\right) \left( \sum_{i\neq n} a_{i\bar i}\right)^2\geq 0, \eeq 
and so
\beq   -\frac{3}{16}|a_{n\bar n}|^2 - \frac{1}{2}\left(\sum_{i\neq n}a_{i\bar i}\right)a_{n\bar n}+ \frac{9}{16(n-1)} \left( \sum_{i\neq n} a_{i\bar i}\right)^2 \geq - \kappa_0\left(\sum_{i=1}^n a_{i\bar i}\right)^2. \eeq 
 Hence,
\beq -\frac{9}{8}\sum_{i=1}^n |a_{n\bar i}|^2 + \frac{7}{8}|a_{n\bar n}|^2 - \frac{1}{2}\sum_{i=1}^n a_{i\bar i}a_{n\bar n} \geq - \kappa_0\left(\sum_{i=1}^n a_{i\bar i}\right)^2- \frac{9}{16}\sum_{i,j=1}^n |a_{i\bar j}|^2. \eeq
This completes the proof of \eqref{claim}.
\eproof

\vskip 1\baselineskip

\section{A Bochner-Kodaira formula for holomorphic sectional curvature}

In this section we derive a Bochner-Kodaira formula tailored to holomorphic sectional curvature, a result of independent significance. The proof strategy underlying Theorem \ref{MainTheoremA} constitutes a logarithmic modification of the  formula \eqref{key7} in Theorem \ref{identity}.

\bproof[Proof of Theorem \ref{identity}] We consider the $(1,1)$-form $ \alpha = \p\omega_2 + \bar\p\omega_1 $. By using similar computations as in the proof of Theorem \ref{MainTheoremA}, we have
\be (\phi\wedge \bar\phi,\alpha) & = & (\alpha,\phi\wedge \bar\phi) = \left(\{R^E\phi,\phi\} - \{\p_E\phi,\p_E\phi\} - \{\bar\p_E\phi,\bar\p_E\phi\}, \phi \wedge \bar\phi \right) \\
& = & \int_M -R(V,\bar V,V,\bar V) \frac{\omega^n}{n!} - \left( \{\p_E\phi,\p_E\phi\},\phi\wedge \bar\phi\right) - \left( \{\bar\p_E\phi,\bar\p_E\phi\}, \phi\wedge \bar\phi\right) \\
& = & \int_M -R(V,\bar V,V,\bar V) \frac{\omega^n}{n!} - (\omega_1,\omega_1)+ (\omega_2,\omega_2), \ee
where the third identity follows from Lemma \ref{omegatwo}. On the other hand, by Lemma \ref{dstarcompute},
\be (\phi\wedge \bar\phi,\alpha) & = & \left( \phi \wedge \bar\phi, \p\omega_2 + \bar\p\omega_1 \right) =
\big( \p^* \left(\phi \wedge \bar\phi\right), \omega_2\big) + \left( \bar\p{}^* \left(\phi \wedge \bar\phi\right), \omega_1\right) \\
& = & (\Delta_{\bp}f\cdot \bar \phi, \omega_2) - (\bar\omega_1,\omega_2) + (-\Delta_{\bp} f\cdot \phi,\omega_1) +(\omega_1,\omega_1).
\ee
Hence, we get the identity
\beq \label{equationA}(\omega_2+\bar{\omega_1},\omega_2) - (\Delta_{\bar\p}f \cdot \bar\phi,\omega_2) = \int_M R(V,\bar V,V,\bar V) \frac{\omega^n}{n!} - (\Delta_{\bar\p}f \cdot \phi,\omega_1) + 2(\omega_1,\omega_1). \eeq
Moreover, a straightforward computation shows
\be
(\Delta_{\bar\p}f \cdot \phi, \p|\phi|^2) & = & \lambda\left(\p^*(f\phi),|\phi|^2\right) = \lambda\left(\smo \Lambda \bar\p(f\phi),|\phi|^2\right) \\
& = & \smo\lambda\left( \Lambda (\bar\phi \wedge \phi + f\bar\p\phi),|\phi|^2\right) \\
& = & -\lambda\int_M |\phi|^4 \frac{\omega^n}{n!} + \lambda^2 \int_M f^2|\phi|^2 \frac{\omega^n}{n!}.
\ee
Hence,
\beq \label{equationB} \lambda\int_M |\phi|^4 \frac{\omega^n}{n!} = (\Delta_{\bar\p}f \cdot \phi, - \omega_1 - \bar\omega_2 + \Delta_{\bar\p}f \cdot \phi), \eeq
On the other hand, we have
\beq  \Delta_{\p}(\bar\p f)=\bp \left(\Delta_{\p} f\right)=\bp(\Delta_{\bp }f)=\lambda \bar \phi.\eeq
Therefore,
\be
\int_M |\phi|^2 |\p\bar\p f|^2 \frac{\omega^n}{n!}  =  (|\phi|^2\p\bar\p f,\p \bar\p f) &=& \left(\p\left(|\phi|^2\bar\phi\right) - \p|\phi|^2 \wedge \bar\phi, \p\bar\p f\right) \\
& = & (|\phi|^2\bar\phi, \Delta_{\p}(\bar\p f)) + (\bar\phi \wedge \p|\phi|^2,\p\bar\p f) \\
& = & \int_M \lambda|\phi|^4 \frac{\omega^n}{n!} + (\p|\phi|^2,I_{\bar V}\p\bar\p f) \\
& = & \int_M \lambda|\phi|^4\frac{\omega^n}{n!} - (\omega_1 + \bar\omega_2, \bar\omega_2).
\ee
Hence,
\beq \label{equationC}\lambda\int_M |\phi|^4\frac{\omega^n}{n!} = \int_M |\phi|^2 |\p\bar\p f|^2 \frac{\omega^n}{n!} + (\omega_1 + \bar\omega_2, \bar\omega_2). \eeq
The summation of equations \eqref{equationA}, \eqref{equationB} and \eqref{equationC} gives
\be 2\lambda \int_M |\phi|^4 \frac{\omega^n}{n!} & = & \int_M \left( R(V,\bar V,V,\bar V) + |\phi|^2 |\p\bar\p f|^2\right) \frac{\omega^n}{n!} + (\omega_1-\Delta_{\bar\p}f\cdot\phi,\omega_1-\Delta_{\bar\p}f\cdot\phi)\\
& &+(\omega_1,\omega_1) + \left[(\Delta_{\bar\p}f \cdot \bar\phi,\omega_2) - (\Delta_{\bar\p}f \cdot \phi,\bar\omega_2)\right] \\
& & +\left[(\omega_1,\bar\omega_2) - (\bar\omega_1,\omega_2)\right] + \left[ (\omega_1,\Delta_{\bar\p}f \cdot \phi) - (\Delta_{\bar\p}f \cdot \phi,\omega_1)\right].\ee
By taking the real part of this equation, 
$$ \left((\Delta_{\bar\p}f \cdot \bar\phi,\omega_2) - (\Delta_{\bar\p}f \cdot \phi,\bar\omega_2)\right) +\left((\omega_1,\bar\omega_2) - (\bar\omega_1,\omega_2)\right) + \left( (\omega_1,\Delta_{\bar\p}f \cdot \phi) - (\Delta_{\bar\p}f \cdot \phi,\omega_1)\right)=0,$$
and we obtain the desired identity \eqref{key7}. Since $ \phi =\omega_g(\cdot, \bar V)$, the inequality \eqref{inequality} is obvious.
\eproof

\noindent A straightforward application of Theorem \ref{identity} is that:
\bcorollary Let $ (M,\omega_g) $ be a complete K\"ahler manifold.
If  its holomorphic sectional curvature satisfies $\mathrm{HSC}  \geq 2 $, then the first
eigenvalue  of the Laplacian satisfies  $ \lambda_1 \geq 2.$\ecorollary

\vskip 2\baselineskip


\end{document}